\documentclass[12pt]{amsart}
\usepackage{amsmath,amssymb}
\usepackage{graphicx}
\usepackage{psfrag}
\usepackage{epsfig}

\begin{document}

\begin{center}
{\large\bf Cutsets in ${\mathcal P}(X)$}
\end{center}

\bigskip
\begin{center}
John Ginsburg\footnote{Winnipeg, Manitoba, Canada} and Bill Sands\footnote{Qualicum Beach, British Columbia, Canada}
\end{center}

\bigskip
\date{\today}

\begin{center}
{\bf Abstract}
\end{center}

\medskip
{\small For any set $X$, ${\mathcal P}(X)$ denotes the collection of all subsets of $X$, ordered by inclusion. A {\it cutset} in ${\mathcal P}(X)$ is a subset of ${\mathcal P}(X)$ which meets every maximal chain of ${\mathcal P}(X)$. A cutset is non-trivial if it does not contain $X$ or the empty set. Our main result is the following. 

\medskip
{\bf Theorem 1}:  Let $X$ be an infinite set of cardinality $\kappa$. Every non-trivial 
cutset in ${\mathcal P}(X)$ contains a chain of cardinality $\kappa^+$ and an antichain of 
cardinality $2^{\kappa}$. }

\vspace{1cm}
\section {\bf Introduction} For any set $X$, ${\mathcal P}(X)$ denotes the collection of all subsets of $X$, ordered by inclusion. A {\it cutset} in ${\mathcal P}(X)$ is a subset of ${\mathcal P}(X)$ which meets every maximal chain of ${\mathcal P}(X)$. A cutset is non-trivial if it does not contain $X$ or the empty set. 

\medskip
Since a maximal chain in a partially ordered set (poset) extends from the bottom to the top,  intuitively a cutset should extend across the poset, in the same way a maximal antichain does: in particular, one might wonder briefly whether every cutset must contain a maximal antichain. Although it is easy to find examples of posets where this fails, in [2] it was proved that, if $X$ is finite, every cutset of ${\mathcal P}(X)$ does contain a maximal antichain. 

\medskip
In this paper we prove a result along similar lines for infinite sets $X$. Our main result is the following. 

\medskip
{\bf Theorem 1}:  Let $X$ be an infinite set of cardinality $\kappa$. Every non-trivial 
cutset in ${\mathcal P}(X)$ contains a chain of cardinality $\kappa^+$ and an antichain of 
cardinality $2^{\kappa}$. 

\medskip
The second part of this theorem --- that non-trivial cutsets in ${\mathcal P}(X)$ must contain large antichains --- supports the image of a cutset spreading ``horizontally" across a poset. But the first part ---  that non-trivial cutsets in ${\mathcal P}(X)$ must contain large chains --- came as somewhat of a surprise.

\medskip
Our notation and terminology involving cardinal and ordinal numbers are standard. By definition, ordinal numbers are transitive sets that are well-ordered by the membership relation. Recall that for any ordinals $\alpha,\beta$, $\alpha<\beta$ and $\alpha\in\beta$ mean the same, and for any ordinal $\beta$, $\beta=\{\alpha|\alpha\mbox{ is an ordinal and } \alpha<\beta\}$. 
Cardinal numbers are initial ordinal numbers --- ordinal numbers that are not equipotent to any smaller ordinal numbers.  For any infinite cardinal number $\kappa$, $\kappa^+$ denotes the first cardinal which is larger than $\kappa$. The reader is referred to [5] for basic set theoretic concepts.

\bigskip
\section{\bf The proof of Theorem 1, with lemmas} The proof is accomplished via three lemmas, as follows.

\medskip 
{\bf Lemma 1}: Let $\kappa$ be an infinite cardinal and let $A$ be a subset of $\kappa$. Let 
$${\mathcal C}_A=\{\emptyset,\kappa\}\cup\{A-\alpha:\alpha<\kappa\}\cup\{ A\cup\alpha:\alpha<\kappa\}.$$
Then ${\mathcal C}_A$  is a maximal chain in ${\mathcal P}(\kappa)$ containing $A$.  

\medskip
{\it Proof}. 
Let $S$ be a subset of $\kappa$ such that $S$ is comparable under inclusion to every member of ${\mathcal C}_A$. We claim that $S\in{\mathcal C}_A$. We assume that $A\subset S$; the argument for $S\subset A$ is dual. Let $\beta$ be the first element of $\kappa$ which is not in $S$. (Otherwise, if there is no such $\beta$, $S=\kappa$ and we are done.) So $\beta\not\in S$ (and therefore $\beta\not\in A$). We claim that $S=A\cup\beta$ (and so $S\in{\mathcal C}_A$).  First, $A\cup\beta\subseteq S$ because $A\subseteq S$ and $\beta$ is the smallest element not in $S$.   So we need to show $S\subseteq A\cup\beta$. If not, let $\gamma\in S-(A\cup\beta)$. Then $\gamma\ge\beta$, so $\gamma>\beta$.   Now consider the set $A\cup\gamma$, which belongs to ${\mathcal C}_A$. Note that $S\not\subseteq A\cup\gamma$, because $\gamma\in S$. Also $A\cup\gamma\not\subseteq S$, because $\beta\in A\cup\gamma$. But then $S$ is not comparable to $A\cup\gamma$, contrary to $S$ being comparable to all members of ${\mathcal C}_A$. \qquad$\Box$

\medskip
{\bf Lemma 2}: Let $\kappa$ be an infinite cardinal and let $\mathcal C$  be a non-trivial cutset in ${\mathcal P}(\kappa)$. For every subset $A$ of $\kappa$, either there is an $\alpha < \kappa$ for which $A - \alpha \in {\mathcal C}$  or there is an $\alpha < \kappa$ for which $A \cup \alpha \in {\mathcal C}$. 

\medskip
{\it Proof}. Follows directly from Lemma 1.\qquad$\Box$

\medskip
{\bf Lemma 3}: Let $\kappa$ be an infinite cardinal. Then 

\smallskip
(a) ${\mathcal P}(\kappa)$ contains a chain ${\mathcal A}$  of cardinality $\kappa^+$  such that $|A-B| = \kappa$ for any two members $A$ and $B$ of ${\mathcal A}$  with $B \subsetneq A$;

\smallskip
(b) ${\mathcal P}(\kappa)$ contains an antichain ${\mathcal B}$ of cardinality $2^{\kappa}$ such that  $|A-B| = \kappa$ for any two distinct members $A$ and $B$ of ${\mathcal B}$. 

\medskip
{\it Proof}. (a) Let $\lambda$ be the smallest cardinal for which $\kappa^{\lambda}>\kappa$. Let $T=\left(\bigcup_{\alpha<\lambda}\kappa^{\alpha}\right)\cup\kappa^{\lambda}$ be the tree of functions, ordered by extension. (Recall that, for cardinals $\alpha$ and $\beta$, the symbol $\alpha^{\beta}$ has two meanings; it is a cardinal itself, but also, as in the definition of $T$, $\alpha^{\beta}$ denotes the set $\{f|f:\beta\to \alpha\}$.)

\medskip
Let $D=\bigcup_{\alpha<\lambda}\kappa^{\alpha}$. Note that $|D|=\kappa$ by the definition of $\lambda$. Define a linear ordering $\le$ on $T$ as follows: for $f,g\in T$, $f\le g$ if either $g$ is an extension of $f$, or if $f(i)<g(i)$ where $i$ is the smallest ordinal for which $f(i)\ne g(i)$. Then

\smallskip\noindent
(i) $\le$ is a linear order on $T$;

\noindent
(ii) $D$ is dense in $T$ (in the order topology);

\noindent
(iii) for any $f_1,f_2\in \kappa^{\lambda}$ with $f_1<f_2$, $|\{d\in D:f_1<d<f_2\}|=\kappa$. 

\smallskip
(i) is easy to prove, and (ii) follows from (iii). For (iii),  let $f$ and $g$ be elements of $\kappa^{\lambda}$ with $f < g$. Let $i_0$ be the smallest element of  $\lambda$ with $f(i_0)\ne g(i_0)$. Then $f(i_0) < g(i_0)$. Let $i_1$ be the next ordinal after $i_0$, and let $y_1$ be any element of $\lambda$ which is larger than $f(i_1)$.
Now, let $\alpha$ be any ordinal with $i_1 < \alpha < \lambda$ and let $h$ be any function from $\alpha$ to $\kappa$ which agrees with $f$ on $i_0 \cup \{i_0\}$ and has the value $y_1$ at $i_1$. In particular, $h(i_0) = f(i_0)$ so $h < g$. Also, the first ordinal $i$ at which $f$ and $h$ have a different value is $i = i_1$. Since $f(i_1) < h(i_1) = y_1$, this implies that $f < h$. So $f < h < g$. For any $\alpha$ for which $i_1 < \alpha < \lambda$, such functions $h$ can have values $h(i)$ defined arbitrarily for ordinals $i > i_1$. So there will be at least $\kappa$
such functions $h$.

\smallskip
Now let $E$ be any subset of $\kappa^{\lambda}$ of cardinality $|E|=\kappa^+$. For each $f\in E$, let $C_f=\{d\in D|d<f\}$. Then ${\mathcal A}=\{C_f|f\in E\}$ is a chain of subsets of $D$ with the properties we want. Since $|D|=\kappa$, we can identify $D$ as $\kappa$ (up to a bijection) and this proves part (a).

\medskip
(b) $X$ can be written as a union of $\kappa$ disjoint sets $A_{\alpha},\alpha<\kappa$ with $|A_{\alpha}|=\kappa$ for all $\alpha$. And each set $A_{\alpha}$ can be written as a union of $\kappa$ disjoint sets $A_{\alpha,\beta}$ with $\beta<\kappa$ such that $|A_{\alpha,\beta}|=\kappa$ for all $\beta$. So $X=\bigcup_{\alpha<\kappa}A_{\alpha}$ and $A_{\alpha}=\bigcup_{\beta<\kappa}A_{\alpha,\beta}$. Now, if $f$ is any function from $\kappa$ to $\kappa$, define the set $S_f$ as follows: $S_f=\bigcup_{\alpha<\kappa}A_{\alpha,f(\alpha)}$. (We take the $f(\alpha)^{\mbox{\tiny th}}$ set from $\{A_{\alpha,\beta}:\beta<\kappa\}$.) Then $\{S_f|f:\kappa\to\kappa\}$ is the desired antichain.\qquad$\Box$

\bigskip
{\it Remark}. The proof of the existence of chains of the desired size in Lemma 3(a) uses a tree of ordinal functions to construct a large
linear order with a small dense subset. This is a standard construction in set theory and apparently originates in [4]. It is worth noting that the connection between chains in ${\mathcal P}(X)$ and dense subsets of linear orderings also extends in the opposite direction. The following are equivalent for infinite cardinals
$\kappa$ and $\lambda$ with $\kappa < \lambda \le 2^{\kappa}$:

\smallskip
(i) If $X$ is a set of cardinality $\kappa$ then ${\mathcal P}(X)$ contains a chain of cardinality $\ge \lambda$;

\smallskip
(ii) there is a densely ordered linearly ordered set of size $\ge \lambda$ which has a dense subset of size $\kappa$

\medskip
 {\it Proof}. (ii)$\to$(i): This is contained in the proof of Lemma 3(a). 

\smallskip
(i)$\to$ (ii): Start with a chain in ${\mathcal P}(X)$ of size $\ge \lambda$. 
We extend it to a maximal chain $M$ in ${\mathcal P}(X)$. Then $M$ under inclusion is the linear order we want in (ii). It has a dense subset of size $\kappa = |X|$ consisting
of the members $M(x)$ of $M$, for $x\in X$, where $M(x)$ is the smallest set in $M$ which contains $x$. ($M(x)$ is the intersection of all members of $M$ which contain $x$. $M$ is closed under intersections because of maximality.)\qquad$\Box$

\bigskip
{\it Proof of Theorem 1}. For the first part, by Lemma 3(a) there is a chain ${\mathcal A}$ of subsets of $\kappa$ such that $|{\mathcal A}|={\kappa}^+$ and $|A-B|=\kappa$ for any sets $A,B\in{\mathcal A}$ for which $A\supset B$. By Lemma 2, for each $A \in {\mathcal A}$ either there is an ordinal $\alpha$ such that $A\cup\alpha\in C$ or there is an ordinal $\alpha$ such that $A-\alpha\in C$. Let ${\mathcal A}_1=\{A\in{\mathcal A}|A\cup\alpha\in{\mathcal C}\mbox{ for some }\alpha<\kappa\}$ and let ${\mathcal A}_2=\{A\in{\mathcal A}|A-\alpha\in{\mathcal C}\mbox{ for some }\alpha<\kappa\}$. Then ${\mathcal A}={\mathcal A}_1\cup{\mathcal A}_2$. At least one of ${\mathcal A}_1,{\mathcal A}_2$ has cardinality $\kappa^+$. We do the case when $|{\mathcal A}_1|=\kappa^+$. So, for every $A\in{\mathcal A}_1$ there is an $\alpha$ for which $A\cup\alpha\in{\mathcal C}$. There are $\kappa^+$ sets $A\in{\mathcal A}_1$ and $\kappa$ elements $\alpha<\kappa$, so by the pigeonhole principle there is an $\alpha_0<\kappa$ and a subcollection ${\mathcal A}'_1$ of ${\mathcal A}_1$ for which $A\cup\alpha_0\in{\mathcal C}$ for all $A\in{\mathcal A}'_1$ and $|{\mathcal A}'_1|=\kappa^+$. ${\mathcal A}'_1\subseteq{\mathcal A}$, so ${\mathcal A}'_1$ is a chain under inclusion. Let $A_1,A_2$ be any two different members of ${\mathcal A}'_1$ with $A_1\subset A_2$. Then $A_1\cup\alpha_0\subseteq A_2\cup\alpha_0$. We cannot have $A_1\cup\alpha_0=A_2\cup\alpha_0$ because this implies that $A_2-A_1\subseteq\alpha_0$ and hence  $|A_2-A_1|\le|\alpha_0|<\kappa$. But $|A_2-A_1|=\kappa$. So $A_1\cup\alpha_0\ne A_2\cup\alpha_0$, so $\{A\cup\alpha_0|A\in{\mathcal A}'_1\}$ is a chain in ${\mathcal C}$ of size $\kappa^+$. 

\medskip
For the second part, by Lemma 3(b) there is an antichain ${\mathcal A}$ of ${\mathcal P}(\kappa)$of cardinality $2^{\kappa}$. Again using Lemma 2,  for each $A\in{\mathcal A}$ either there is an ordinal $\alpha$ such that $A\cup\alpha\in{\mathcal C}$ or there is an ordinal $\alpha$ such that $A-\alpha\in{\mathcal C}$. We again let
${\mathcal A}_1 = \{A\in{\mathcal A}|A\cup\alpha\in{\mathcal C}\mbox{ for some } \alpha < \kappa\}$ and ${\mathcal A}_2 = \{A\in{\mathcal A}|A-\alpha\in{\mathcal C}\mbox{ for some } \alpha < \kappa\}$. As above, ${\mathcal A} = {\mathcal A}_1 \cup {\mathcal A}_2$ and at least one of ${\mathcal A}_1, {\mathcal A}_2$ has cardinality $2^{\kappa}$. We again do the case when it is ${\mathcal A}_1$.
 So, for each $A$ in ${\mathcal A}_1$, there is an $\alpha < \kappa$ with $A \cup \alpha \in {\mathcal C}$. For each $\alpha < \kappa$, let ${\mathcal B}_{\alpha} = \{A \in {\mathcal A}_1 | A \cup \alpha \in {\mathcal C}\}$. 
So we have ${\mathcal A}_1 = \bigcup_{\alpha<\kappa}{\mathcal B}_{\alpha}$. 
But ${\mathcal A}_1$ has cardinality $2^{\kappa}$ and the cofinality of $2^{\kappa}$ is larger than $\kappa$, by K\"{o}nig's theorem. (See [5], Theorem 14 in Chapter VI, for a proof of K\"{o}nig's theorem.)  So at least one of the sets ${\mathcal B}_{\alpha}$ must have cardinality $2^{\kappa}$, say
${\mathcal B}_{\alpha_0}$. So we have $2^{\kappa}$ sets $A \in {\mathcal A}$ for which $A\cup \alpha_0$
is in ${\mathcal C}$. Since $|A_1 - A_2| = \kappa$ for any two sets in ${\mathcal A}$, it follows that no two sets $A_1 \cup \alpha_0, A_2 \cup \alpha_0$ for $A_1$ and $A_2$ in ${\mathcal A}$, can be comparable under inclusion: if $A_1 \cup \alpha_0$ were a subset of $A_2 \cup \alpha_0$, then we would have $A_1 - A_2$ a subset of $\alpha_0$, which is impossible since $|A_1 - A_2| = \kappa$. So these $2^{\kappa}$ sets form an antichain which is contained in the cutset ${\mathcal C}$. \qquad$\Box$

\medskip
{\it Note}: We note that, because of related independence results, for example see [3], Section 4, we cannot replace $\kappa^+$ by $2^{\kappa}$ in Lemma 3(a) or in the first conclusion of Theorem 1.

\bigskip
\section{\bf Further results and questions} 

\medskip
We first remark that, in the countable case $\kappa = \omega$, Theorem 1 can be sharpened. 

\medskip
{\bf Theorem 2}. In ${\mathcal P}(\omega)$, every non-trivial cutset must contain a chain of 
cardinality $2^{\omega}$ as well as an antichain of cardinality $2^{\omega}$. 

\medskip
{\it Proof}. From the proof of Theorem 1 and our remarks following it on dense subsets of linear orders, it is enough to note that there is a linear order of cardinality $2^{\omega}$ (the real numbers) which has a dense subset of cardinality $\omega$ (the rationals). \qquad$\Box$

\medskip
Theorem 1 shows that, for infinite sets $X$, all non-trivial cutsets in ${\mathcal P}(X)$ must contain large chains. This is in stark contrast to the case when $X$ is finite, where ${\mathcal P}(X)$ contains non-trivial antichain cutsets --- its levels. 

\medskip
A cutset ${\mathcal C}$ is {\it minimal} if no proper subset of ${\mathcal C}$ is a cutset. Via the levels, every element of ${\mathcal P}(X)$ is in a minimal cutset of ${\mathcal P}(X)$ when $X$ is finite. We next prove that this property also holds in the infinite case.

\medskip
{\bf Theorem 3}. Let $X$ be infinite. Then any set $A\in{\mathcal P}(X)$ belongs to a minimal cutset of ${\mathcal P}(X)$. 

\medskip
{\it Proof}. Let $x,y$ be distinct elements of $X$, and let  
$${\mathcal C}=\{S\in{\mathcal P}(X):|S\cap\{x,y\}|=1\}.$$ 
First we show ${\mathcal C}$ is a cutset in ${\mathcal P}(X)$; in other words, ${\mathcal C}$ meets every maximal chain of ${\mathcal P}(X)$. Let $M$ be any maximal chain in ${\mathcal P}(X)$. We want to show that ${\mathcal C}$ contains an element of $M$. Let $A$ be the intersection of all of the sets belonging to $M$ which contain both $x$ and $y$. ($X$ is one such set.) We note that $A$ is an element of $M$, by maximality. In fact, because $M$ is maximal, it is closed under arbitrary unions and intersections. Now let $U$ be the union of all the members $B$ of $M$ for which $B$ is a proper subset of $A$. Therefore $U$ is a subset of $A$. We cannot have both $x$ and $y$ in $U$, because there would then be members $B_1$ and $B_2$ of $M$, both proper subsets of $A$, with $x\in B_1$ and $y\in B_2$. One of the sets $B_1$ or $B_2$ must contain the other. Whichever one it is, it would be a strictly smaller member of $M$ containing both $x$ and $y$, which is impossible. So let us assume that $x$ is not in $U$. In particular, $U$ is a proper subset of $A$. So $U$ must be the immediate predecessor of $A$ in $M$. Since $U\subset A -\{x\}\subset A$, we must have $U = A - \{x\}$, since otherwise we could add the set $A - \{x\}$ to $M$ to form a larger chain of sets, contradicting maximality. So $U$ is a member of $M$ which contains $y$ but not $x$.  Hence $U$ is an element of $M$ which belongs to ${\mathcal C}$.

\medskip
Now we show that ${\mathcal C}$ is minimal. Let $S_0$ be any element of ${\mathcal C}$. We can assume that $x\in S_0$ and $y\not\in S_0$. We want to show that ${\mathcal C} - \{S_0\}$ is not a cutset in ${\mathcal P}(X)$.  Consider the three-element chain  $G = \{ S_0 - \{x\}, S_0, S_0 \cup \{y\}\}$. There is a maximal chain $M$ in ${\mathcal P}(X)$ with $G$ contained in $M$. Any set which belongs to $M$ and which is strictly larger than $S_0$ contains both $x$ and $y$, and hence does not belong to ${\mathcal C}$. Any set which belongs to $M$ and which is strictly smaller than $S_0$  must be contained in $S_0 - \{x\}$, so it doesn't contain $x$, nor does it contain $y$, because $S_0$ doesn't contain $y$. So it doesn't belong to ${\mathcal C}$ either. Thus the only member of $M$ which does belong to ${\mathcal C}$ is $S_0$. Therefore ${\mathcal C} - \{S_0\}$ does not meet $M$. Therefore ${\mathcal C} - \{S_0\}$ is not a cutset.

\medskip
Now let $S$ be any subset of $X$. If  $S = X$ or $S = \emptyset$, the conclusion
is obvious. Otherwise, there must be at least one element $x\in S$ and there must be at least one element $y\in X - S$.  For this pair of elements $x$ and $y$, the cutset ${\mathcal C}$ defined above is a minimal cutset containing $S$.\qquad$\Box$

\bigskip
Theorem 1 shows that, for infinite $X$, non-trivial cutsets of ${\mathcal P}(X)$ must contain ``large" antichains, whereas [2] has a proof that, for finite $X$, cutsets of ${\mathcal P}(X)$ must contain maximal antichains.  We have been unable to resolve the following, even for the countable case:

\medskip
{\it Question}. For $X$ infinite, does every cutset of ${\mathcal P}(X)$ contain a maximal antichain of ${\mathcal P}(X)$? 

\medskip
To end, a personal note: We are two former university professors, from the University of Winnipeg and the University of Calgary respectively, both retired for over ten years. This is
our seventh joint paper together, but our first since retirement. Our initial interest in cutsets began long ago, originating from [1]. That paper was concerned with infinite  sets, but it has led both of us, and others, to many related investigations concerning finite partially ordered sets and cutsets. After many years, this paper returns us to the study of cutsets in the realm of infinite sets.

\medskip
  Reactions to this paper are very much encouraged! Please direct  comments (criticisms, compliments, corrections, whatever!) to either of us at ginsburg.john@gmail.com or sands@ucalgary.ca.

\bigskip
{\bf References:}

\medskip
[1] M.~Bell and J.~Ginsburg, Compact spaces and spaces of maximal complete subgraphs, {\it Trans. Amer. Math. Soc.} {\bf 283} (1984), 329--339. 

\smallskip
[2] D.~Duffus, B.~Sands and P.~Winkler, Maximal chains and antichains in Boolean lattices, {\it SIAM Journal on Discrete Mathematics} {\bf 3}(2) (1990), 197--205.

\smallskip
[3] W.~Mitchell, Aronszajn trees and the independence of the transfer property, {\it Annals Of Mathematical Logic} {\bf 5} (1972),  21--46. 

\smallskip
[4] Ð.~Kurepa. Ensembles ordonnés et ramifiés, {\it Publ. Math. Univ. Belgrade} {\bf 4} (1935), 1--138.

\smallskip
[5] W.~Sierpi\'{n}ski, {\it Cardinal and Ordinal Numbers} (1965), PWN - Polish Scientific Publishers,  in English - 2nd ed. revised. 

\end{document}